# Fractional Weierstrass function by application of Jumarie fractional trigonometric functions and its analysis

Uttam Ghosh[1], Susmita Sarkar[2] and Shantanu Das[3]

[1] Department of Mathematics, Nabadwip Vidyasagar College, Nabadwip, Nadia, West Bengal, India;
email : uttam_math@yahoo.co.in

[2] Department of Applied Mathematics, University of Calcutta, Kolkata, India
email : susmita62@yahoo.co.in

[3] Scientist H+, Reactor Control Systems Design Section E & I Group BARC Mumbai India
Senior Research Professor, Dept. of Physics, Jadavpur University Kolkata
Adjunct Professor. DIAT-Pune
UGC Visiting Fellow. Dept of Appl. Mathematics; Univ. of Calcutta
email : shantanu@barc.gov.in

## Abstract

The classical example of no-where differentiable but everywhere continuous function is Weierstrass function. In this paper we define the fractional order Weierstrass function in terms of Jumarie fractional trigonometric functions. The Holder exponent and Box dimension of this function are calculated here. It is established that the Holder exponent and Box dimension of this fractional order Weierstrass function are the same as in the original Weierstrass function, independent of incorporating the fractional trigonometric function. This is new development in generalizing the classical Weierstrass function by usage of fractional trigonometric function and obtain its character and also of fractional derivative of fractional Weierstrass function by Jumarie fractional derivative, and establishing that roughness index are invariant to this generalization.

## 1.0 Introduction

Fractional geometry, fractional dimension is an important branch of science to study the irregularity of a function, graph or signals [1-3]. On the other hand fractional calculus is another developing mathematical tool to study the continuous but non-differentiable functions (signals) where the conventional calculus fails [4-11]. Many authors are trying to relate between the fractional derivative and fractional dimension [1, 12-15]. The functions which are continuous but non-differentiable in integer order calculus can be characterized in terms of fractional calculus and especially through Holder exponent [10, 16]. To study the no-where differentiable functions authors in [12-16] used different type of fractional derivatives. Jumarie [17] defines the fractional trigonometric functions in terms of Mittag-Leffler function and established different useful fractional trigonometric formulas. The fractional order derivatives of those functions were established in-terms the Jumarie [17-18] modified fractional order derivatives. In this paper we define the fractional order Weierstrass functions in terms the fractional order sine function. The Holder exponent, box-dimension (Fractional dimension) of graph of this function is obtained here; also the fractional order derivative of this function is established here. This is new development in generalizing the classical Weierstrass function by usage of fractional trigonometric function and obtain its character. The paper is organized as sections; with section-2 dealing with describing Jumarie fractional derivative and Mittag-Leffler function of one and two parameter type, fractional trigonometric function of one and two parameter type and their Jumarie fractional derivatives are derived. In this section we derived useful relations of fractional trigonometric function that we shall be using for our calculations-in characterizing fractional Weierstrass function. We continue this section by introducing Lipschitz Holder exponent (LHE)- its definition, its relation to Hurst exponent and fractional dimension and definition of Holder continuity, and we define here the classical Weierstrass function. These parameters are basic parameter to indicate roughness index of a function or graph. In section-3 we describe the fractional Weierstrass function by generalizing the classical Weierstrass function by use of fractional sine trigonometric function. Subsequently we apply derived identities of fractional trigonometric functions to evaluate the properties of this new fractional Weierstrass function. In section-4 we do derivation of properties of fractional derivative of fractional Weierstrass function, and conclude the paper with conclusion and references.

## 2.0 Jumarie fractional order derivative and Mittag-Leffler Function

## a) Fractional order derivative of Jumarie Type

Jumarie [17] defined the fractional order derivative modifying the Left Riemann-Liouvellie (RL) fractional derivative in the form for the function $f(x)$ in the interval $a$ to $x$, with $f(x) = 0$ for $x < a$.

$$
{}_a^J D_x^\alpha [f(x)] = \begin{cases} \dfrac{1}{\Gamma(-\alpha)} \displaystyle\int_a^x (x-\tau)^{-\alpha-1} f(\tau) d\tau, & \alpha < 0. \\ \dfrac{1}{\Gamma(1-\alpha)} \dfrac{d}{dx} \displaystyle\int_a^x (x-\tau)^{-\alpha} [f(\tau) - f(a)] d\tau, & 0 < \alpha < 1 \\ \left( f^{(\alpha-m)}(x) \right)^{(m)}, & m \leq \alpha < m+1. \end{cases}
$$

In the above definition, the first expression is just Riemann-Liouvelli fractional integration; the second line is Riemann-Liouvelli fractional derivative of order $0 < \alpha < 1$ of offset function that is $f(x) - f(a)$. For $\alpha > 1$, we use the third line; that is first we differentiate the offset function with order $0 < (\alpha - m) < 1$, by the formula of second line, and then apply whole $m$ order differentiation to it. Here we chose integer $m$, just less than the real number $\alpha$; that is $m \leq \alpha < m+1$. In this paper we use symbol ${}_0^J D_x^\alpha$ to denote Jumarie fractional derivative operator, as defined above. In case the start point value $f(a)$ is un-defined there we take finite part of the offset function as $f(x) - f(a^+)$; for calculations. Note in the above Jumarie definition ${}_0 D_x^\alpha [C] = 0$, where $C$ is constant function, otherwise in RL sense, the fractional derivative of a constant function is ${}_0 D_x^\alpha [C] = C \frac{x^{-\alpha}}{\Gamma(1-\alpha)}$, that is a decaying power-law function. Also we purposely state that $f(x) = 0$ for $x < 0$ in order to have initialization function in case of fractional differ-integration to be zero, else results are difficult [9].

## b) Mittag-Leffler Function

The Mittag-Leffler function [19- 22] of one parameter is denoted by $E_\alpha(x)$ and defined by

$$E_\alpha(x) = \sum_{k=0}^{\infty} \frac{x^k}{\Gamma(1+k\alpha)} \qquad (1)$$

This function plays a crucial role in classical calculus for $\alpha = 1$, for $\alpha = 1$ it becomes the exponential function, that is $\exp(x) = E_1(x)$

$$\exp(x) = \sum_{k=0}^{\infty} \frac{x^k}{k!} \qquad (2)$$

Like the exponential function; $E_\alpha(x^\alpha)$ play important role in fractional calculus. The function $E_\alpha(x^\alpha)$ is a fundamental solution of the Jumarie type fractional differential equation $_0D_y^\alpha[y] = y$, where $_0D_x^\alpha$ is Jumarie derivative operator as described above. The other important function is the two parameter Mittag-Leffler function is denoted $E_{\alpha,\beta}(x)$ and defined by following series

$$E_{\alpha,\beta}(x) = \sum_{k=0}^{\infty} \frac{x^k}{\Gamma(\beta+k\alpha)} \qquad (3)$$

The functions (1) and (3) play important role in fractional calculus, also we note that $E_{\alpha,1}(x) = E_\alpha(x)$. Again from Jumarie definition of fractional derivative we have $_0^JD_x^\alpha[1] = 0$ and $_0^JD_x^\alpha[x^\beta] = \frac{\Gamma(1+\beta)}{\Gamma(1+\beta-\alpha)} x^{\beta-\alpha}$. We now consider the Mittag-Leffler function in the following form in infinite series representation for $f(x) = E_\alpha(x^\alpha)$ for $x \geq 0$ and $f(x) = 0$ for $x < 0$ as;

$$E_\alpha(x^\alpha) = 1 + \frac{x^\alpha}{\Gamma(1+\alpha)} + \frac{x^{2\alpha}}{\Gamma(1+2\alpha)} + \frac{x^{3\alpha}}{\Gamma(1+3\alpha)} + ...$$

Then taking Jumarie fractional derivative of order $0 < \alpha < 1$ term by term for the above series we obtain the following by using the formula $_0^JD_x^\alpha[x^\beta] = \frac{\Gamma(1+\beta)}{\Gamma(1+\beta-\alpha)} x^{\beta-\alpha}$ and $_0^JD_x^\alpha[1] = 0$

$$_0^JD_x^\alpha[E_\alpha(ax^\alpha)] = _0^JD_x^\alpha\left[1 + \frac{ax^\alpha}{\Gamma(1+\alpha)} + \frac{a^2x^{2\alpha}}{\Gamma(1+2\alpha)} + \frac{a^3x^{3\alpha}}{\Gamma(1+3\alpha)} + ...\right]$$

$$= 0 + a + \frac{a^2x^\alpha}{\Gamma(1+\alpha)} + \frac{a^3x^{2\alpha}}{\Gamma(1+2\alpha)} + \frac{a^4x^{3\alpha}}{\Gamma(1+3\alpha)} + ...$$

$$= aE_\alpha(ax^\alpha)$$

Again we derive Jumarie derivative of order $\beta$ for one parameter Mittag-Leffler function $E_\alpha(x^\alpha)$ and thereby get two parameter Mittag-Leffler function. Here also we use for term by term Jumarie derivative ${}_0^J D_x^\alpha \left[ x^\upsilon \right] = \frac{\Gamma(1+\upsilon)}{\Gamma(1+\upsilon-\alpha)} x^{\upsilon-\alpha}$ and ${}_0^J D_x^\alpha [1] = 0$.

$$\begin{aligned}
{}_0^J D_x^\beta \left[ E_\alpha(x^\alpha) \right] &= {}_0^J D_x^\beta \left[ 1 + \frac{x^\alpha}{\Gamma(1+\alpha)} + \frac{x^{2\alpha}}{\Gamma(1+2\alpha)} + \frac{x^{3\alpha}}{\Gamma(1+3\alpha)} + \ldots \right] \\
&= 0 + \frac{x^{\alpha-\beta}}{\Gamma(1+\alpha-\beta)} + \frac{x^{2\alpha-\beta}}{\Gamma(1+2\alpha-\beta)} + \frac{x^{3\alpha-\beta}}{\Gamma(1+3\alpha-\beta)} + \ldots \\
&= x^{\alpha-\beta} E_{\alpha,\alpha-\beta+1}(x^\alpha)
\end{aligned}$$

Jumarie [18] defined the fractional sine and cosine function in the following form

$$E_\alpha(ix^\alpha) \stackrel{\text{def}}{=} \cos_\alpha(x^\alpha) + i\sin_\alpha(x^\alpha)$$

$$\cos_\alpha(x^\alpha) \stackrel{\text{def}}{=} \sum_{k=0}^{\infty} (-1)^k \frac{x^{2k\alpha}}{\Gamma(1+2\alpha k)}$$

$$\sin_\alpha(x^\alpha) \stackrel{\text{def}}{=} \sum_{k=0}^{\infty} (-1)^k \frac{x^{(2k+1)\alpha}}{\Gamma(1+(1+2k)\alpha)}$$

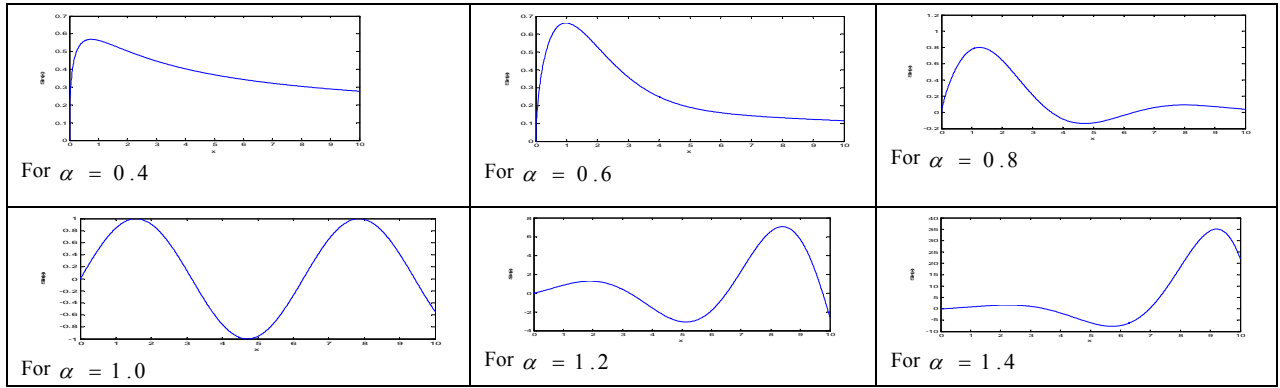

| For $\alpha = 0.4$ | For $\alpha = 0.6$ | For $\alpha = 0.8$ |
| For $\alpha = 1.0$ | For $\alpha = 1.2$ | For $\alpha = 1.4$ |

**Fig.1 Graph of $\sin_\alpha(x^\alpha)$**

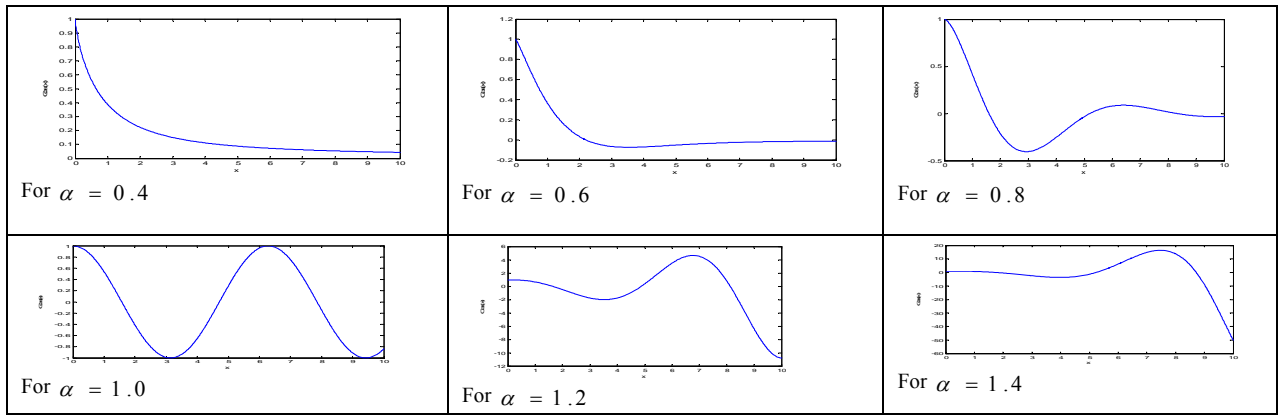

| For $\alpha = 0.4$ | For $\alpha = 0.6$ | For $\alpha = 0.8$ |
| For $\alpha = 1.0$ | For $\alpha = 1.2$ | For $\alpha = 1.4$ |

**Fig.2 Graph of $\cos_\alpha(x^\alpha)$**

From figure-1 and 2 it is observed that for $\alpha < 1$ both the fractional trigonometric functions $\sin_\alpha(x^\alpha)$ and $\cos_\alpha(x^\alpha)$ is decaying functions like damped oscillatory motion. For $\alpha = 1$ it is like simple harmonic motion with sustained oscillations; and for $\alpha > 1$ it grows while it oscillates infinitely; like unstable oscillator.

The series presentation of $f(t) = \cos_\alpha(t^\alpha)$ for $t \geq 0$ and $f(t) = 0$ for $t < 0$ is following

$$\cos_\alpha(at^\alpha) = 1 - \frac{a^2 t^{2\alpha}}{\Gamma(1+2\alpha)} + \frac{a^4 t^{4\alpha}}{\Gamma(1+4\alpha)} - \frac{a^6 t^{6\alpha}}{\Gamma(1+6\alpha)} + \dots$$

Taking term by term Jumarie derivative we get

$$_0^J D_t^\alpha \left[\cos_\alpha(at^\alpha)\right] = 0 - a^2 \frac{\Gamma(1+2\alpha) t^{2\alpha - \alpha}}{\Gamma(1+2\alpha)\Gamma(1+\alpha)} + a^4 \frac{\Gamma(1+4\alpha) t^{4\alpha - \alpha}}{\Gamma(1+4\alpha)\Gamma(1+3\alpha)}$$

$$- a^6 \frac{\Gamma(1+6\alpha) t^{6\alpha - \alpha}}{\Gamma(1+6\alpha)\Gamma(1+5\alpha)} + \dots$$

$$= -a \left[ \frac{at^\alpha}{\Gamma(1+\alpha)} - \frac{a^3 t^{3\alpha}}{\Gamma(1+3\alpha)} + \dots \right]$$

$$= -a \sin_\alpha(at^\alpha)$$

The series presentation of $f(t) = \sin_\alpha(t^\alpha)$, for $t \geq 0$ with $f(t) = 0$ for $t < 0$ is

$$\sin_\alpha(at^\alpha) = \frac{at^\alpha}{\Gamma(1+\alpha)} - \frac{a^3 t^{3\alpha}}{\Gamma(1+3\alpha)} + \frac{a^5 t^{5\alpha}}{\Gamma(1+5\alpha)} - \frac{a^7 t^{7\alpha}}{\Gamma(1+7\alpha)} + \dots$$

Taking term by term Jumarie derivative we get

$$_0^J D_t^\alpha \left[\sin_\alpha(at^\alpha)\right] = a \frac{\Gamma(1+\alpha) t^{\alpha - \alpha}}{\Gamma(1+\alpha)\Gamma(1+\alpha - \alpha)} - a^3 \frac{\Gamma(1+3\alpha) t^{3\alpha - \alpha}}{\Gamma(1+3\alpha)\Gamma(1+3\alpha - \alpha)} + a^5 \frac{\Gamma(1+5\alpha) t^{5\alpha - \alpha}}{\Gamma(1+5\alpha)\Gamma(1+5\alpha - \alpha)}$$

$$- a^7 \frac{\Gamma(1+7\alpha) t^{7\alpha - \alpha}}{\Gamma(1+7\alpha)\Gamma(1+7\alpha - \alpha)} + \dots$$

$$= a \left[ 1 - \frac{a^2 t^{2\alpha}}{\Gamma(1+2\alpha)} + \frac{a^4 t^{4\alpha}}{\Gamma(1+4\alpha)} - \frac{a^6 t^{6\alpha}}{\Gamma(1+6\alpha)} + \dots \right]$$

$$= a \cos_\alpha(t^\alpha)$$

Thus we get

$$_0^J D_x^\alpha \left[\cos_\alpha(ax^\alpha)\right] = -a \sin_\alpha(x^\alpha) \text{ and } _0^J D_x^\alpha \left[\sin_\alpha(ax^\alpha)\right] = a \cos_\alpha(x^\alpha)$$

Jumarie in [18] established $E_\alpha\left(i((x+y)^\alpha)\right) = E_\alpha(ix^\alpha) \times E_\alpha(iy^\alpha)$. Proof of the above relation we reproduce. Let us consider a function $f(x)$ which satisfies the condition

$$f(\lambda x^\alpha) f(\lambda y^\alpha) = f\left(\lambda(x+y)^\alpha\right)$$

Differentiating both side with respect to $x$ and $y$ of $\alpha$-order respectively we get the following

First consider $y$ a constant, and we fractionally differentiate w.r.t. $x$ by Jumarie derivative

$$^J_0D^\alpha_x\left[f(\lambda x^\alpha)\right] \times {^J_0D^\alpha_x}\left[\lambda x^\alpha\right] \times f(\lambda y^\alpha) = {^J_0D^\alpha_x}\left[f\left(\lambda(x+y)^\alpha\right)\right] \times {^J_0D^\alpha_x}\left[\lambda(x+y)^\alpha\right] \times {^J_0D^\alpha_x}[x+y]$$

$$^J_0D^\alpha_x\left[f(\lambda x^\alpha)\right] = f^\alpha(\lambda x^\alpha) \qquad {^J_0D^\alpha_x}\left[f\left(\lambda(x+y)^\alpha\right)\right] = f^\alpha\left(\lambda(x+y)^\alpha\right)$$

$$f^\alpha(\lambda x^\alpha) \times \left[(\lambda)\frac{\Gamma(1+\alpha)x^{\alpha-\alpha}}{\Gamma(1+\alpha-\alpha)}\right] \times f(\lambda y^\alpha) = f^\alpha\left(\lambda(x+y)^\alpha\right) \times \left[(\lambda)\frac{\Gamma(1+\alpha)(x+y)^{\alpha-\alpha}}{\Gamma(1+\alpha-\alpha)}\right] \times {^J_0D^\alpha_x}[x+y]$$

Now we consider $x$ as constant and do the following steps

$$f(\lambda x^\alpha) \times {^J_0D^\alpha_y}\left[f(\lambda y^\alpha)\right] \times {^J_0D^\alpha_y}\left[\lambda y^\alpha\right] = {^J_0D^\alpha_y}\left[f\left(\lambda(x+y)^\alpha\right)\right] \times {^J_0D^\alpha_y}\left[\lambda(x+y)^\alpha\right] \times {^J_0D^\alpha_y}[x+y]$$

$$^J_0D^\alpha_y\left[f(\lambda y^\alpha)\right] = f^\alpha(\lambda y^\alpha) \qquad {^J_0D^\alpha_y}\left[f\left(\lambda(x+y)^\alpha\right)\right] = f^\alpha\left(\lambda(x+y)^\alpha\right)$$

$$f(\lambda x^\alpha) \times f^\alpha(\lambda y^\alpha) \times \left[(\lambda)\frac{\Gamma(1+\alpha)y^{\alpha-\alpha}}{\Gamma(1+\alpha-\alpha)}\right] = f^\alpha\left(\lambda(x+y)^\alpha\right) \times \left[(\lambda)\frac{\Gamma(1+\alpha)(x+y)^{\alpha-\alpha}}{\Gamma(1+\alpha)}\right] \times {^J_0D^\alpha_y}[x+y]$$

Here we put equivalence of $^J_0D^\alpha_y[x+y] \equiv {^J_0D^\alpha_y}[x+y] \equiv {^J_0D^\alpha_u}[u+C]$, with $C$ as constant; that is when $x$ or $y$ are taken as constant the function form of these two quantities gets equivalent that is equivalent to $^J_0D^\alpha_u[u]$ as Jumarie fractional derivative of constant is zero. Therefore the RHS of above two expressions are equal, from that we get the following

$$f^\alpha(\lambda x^\alpha) f(\lambda y^\alpha) = f(\lambda x^\alpha) f^\alpha(\lambda y^\alpha)$$

$$\frac{f^\alpha(\lambda x^\alpha)}{f(\lambda x^\alpha)} = \frac{f^\alpha(\lambda y^\alpha)}{f(\lambda y^\alpha)}$$

The above two may be equated to a constant say $\lambda$. Then we have $f^\alpha(\lambda x^\alpha) = \lambda f(\lambda x^\alpha)$, or we write $^J_0D^\alpha_x\left[f(\lambda x^\alpha)\right] = \lambda f(\lambda x^\alpha)$. From the property of Mittag-Leffler function and Jumarie derivative of the Mittag-Leffler function we know that $^J_0D^\alpha_x\left[E_\alpha(ax^\alpha)\right] = aE_\alpha(ax^\alpha)$; we imply that the solution of $f^\alpha(\lambda x^\alpha) = \lambda f(\lambda x^\alpha)$ is $f(\lambda x^\alpha) = E_\alpha(\lambda x^\alpha)$. Therefore $E_\alpha(\lambda x^\alpha)$ satisfies the condition $f(\lambda x^\alpha) f(\lambda y^\alpha) = f\left(\lambda(x+y)^\alpha\right)$, or $E_\alpha(\lambda x^\alpha) E_\alpha(\lambda y^\alpha) = E_\alpha\left(\lambda(x+y)^\alpha\right)$.
Considering $\lambda = i$, we therefore can write the following identity

$$E_\alpha\left(i((x+y)^\alpha)\right) = E_\alpha\left(i(x^\alpha)\right) E_\alpha\left(i(y^\alpha)\right)$$

Using definition $E_\alpha(ix^\alpha) = \cos_\alpha(x^\alpha) + i\sin(x^\alpha)$ we expand the above as depicted below

$$\cos_\alpha(x+y)^\alpha + i\sin_\alpha(x+y)^\alpha = \left[\cos_\alpha(x^\alpha) + i\sin_\alpha(x^\alpha)\right] \times \left[\cos_\alpha(y^\alpha) + i\sin_\alpha(y^\alpha)\right]$$

$$= \left[\cos_\alpha(x^\alpha)\cos_\alpha(y^\alpha) - \sin_\alpha(y^\alpha)\sin_\alpha(x^\alpha)\right]$$

$$+ i\left[\sin_\alpha(x^\alpha)\cos_\alpha(y^\alpha) + \sin_\alpha(y^\alpha)\cos_\alpha(x^\alpha)\right]$$

Comparing real and imaginary part in above derived relation we get the following

$$\sin_\alpha(x+y)^\alpha = \sin_\alpha(x^\alpha)\cos_\alpha(y^\alpha) + \sin_\alpha(y^\alpha)\cos_\alpha(x^\alpha)$$

$$\cos_\alpha(x+y)^\alpha = \cos_\alpha(x^\alpha)\cos_\alpha(y^\alpha) - \sin_\alpha(y^\alpha)\sin_\alpha(x^\alpha)$$

This is very useful relation as in conjugation with classical trigonometric functions, and we will be using these relations in our analysis of fractional Weierstrass function and its fractional derivative.

Let us define $\cos_{\alpha,\beta}(x^\alpha)$ and $\sin_{\alpha,\beta}(x^\alpha)$ as depicted below

$$\cos_{\alpha,\beta}(x^\alpha) \overset{\text{def}}{=} \sum_{k=0}^{\infty}(-1)^k \frac{x^{2k\alpha}}{\Gamma(\beta+2\alpha k)} = \frac{1}{\Gamma(\beta)} - \frac{x^{2\alpha}}{\Gamma(2\alpha+\beta)} + \frac{x^{4\alpha}}{\Gamma(4\alpha+\beta)} - \ldots$$

$$\sin_{\alpha,\beta}(x^\alpha) \overset{\text{def}}{=} \sum_{k=0}^{\infty}(-1)^k \frac{x^{(2k+1)\alpha}}{\Gamma(\beta+(1+2k)\alpha)} = \frac{x^\alpha}{\Gamma(\alpha+\beta)} - \frac{x^{3\alpha}}{\Gamma(3\alpha+\beta)} + \frac{x^{5\alpha}}{\Gamma(5\alpha+\beta)} + \ldots$$

Now with this and with definition of two parameter Mittag-Leffler (3) function with imaginary argument we get the following useful identity

$$E_{\alpha,\beta}(ix^\alpha) = \sum_{k=0}^{\infty}\frac{(ix^\alpha)^k}{\Gamma(\beta+k\alpha)}$$

$$= \frac{1}{\Gamma(\beta)} + \frac{(ix)^\alpha}{\Gamma(\alpha+\beta)} + \frac{(ix^\alpha)^2}{\Gamma(2\alpha+\beta)} + \frac{(ix^\alpha)^3}{\Gamma(3\alpha+\beta)} + \ldots$$

$$= \left(\frac{1}{\Gamma(\beta)} - \frac{x^{2\alpha}}{\Gamma(2\alpha+\beta)} + \frac{x^{4\alpha}}{\Gamma(4\alpha+\beta)} - \ldots\right) + i\left(\frac{x^\alpha}{\Gamma(\alpha+\beta)} - \frac{x^{3\alpha}}{\Gamma(3\alpha+\beta)} + \frac{x^{5\alpha}}{\Gamma(5\alpha+\beta)} + \ldots\right)$$

$$= \cos_{\alpha,\beta}(x^\alpha) + i\sin_{\alpha,\beta}(x^\alpha)$$

Now for $\beta > \alpha$, we do the Jumarie derivative of order $\alpha$ on the function $f(x) = x^{\beta-1}\cos_{\alpha,\beta}(x^\alpha)$ as depicted in following steps, with formula ${}_0^J D_x^\alpha[x^\upsilon] = \frac{\Gamma(1+\upsilon)}{\Gamma(1+\upsilon-\alpha)}x^{\upsilon-\alpha}$ and ${}_0^J D_x^\alpha[1] = 0$.

$$
\begin{aligned}
{}_0^J D_x^\alpha \left[ x^{\beta-1} \cos_{\alpha,\beta}(x^\alpha) \right] &= {}_0^J D_x^\alpha \left[ (x^{\beta-1}) \times \left( \frac{1}{\Gamma(\beta)} - \frac{x^{2\alpha}}{\Gamma(2\alpha+\beta)} + \frac{x^{4\alpha}}{\Gamma(4\alpha+\beta)} - \ldots \right) \right] \\
&= \frac{x^{\beta-\alpha-1}}{\Gamma(\beta-\alpha)} - \frac{x^{\beta+\alpha-1}}{\Gamma(\beta+\alpha)} + \frac{x^{\beta+3\alpha-1}}{\Gamma(\beta+3\alpha)} - \ldots \\
&= (x^{\beta-\alpha-1}) \times \left[ \frac{1}{\Gamma(\beta-\alpha)} - \frac{x^{2\alpha}}{\Gamma(\beta-\alpha+2\alpha)} + \frac{x^{\beta+3\alpha-1}}{\Gamma(\beta-\alpha+4\alpha)} - \ldots \right] \\
&= x^{\beta-\alpha-1} \cos_{\alpha,\beta-\alpha}(x^\alpha)
\end{aligned}
$$

Thus we get a very useful relation

$$ {}_0^J D_x^\alpha \left[ x^{\beta-1} \cos_{\alpha,\beta}(x^\alpha) \right] = x^{\beta-\alpha-1} \cos_{\alpha,\beta-\alpha}(x^\alpha) $$

Similarly it can be shown that

$$ {}_0^J D_x^\alpha \left[ x^{\beta-1} \sin_{\alpha,\beta}(x^\alpha) \right] = x^{\beta-\alpha-1} \sin_{\alpha,\beta-\alpha}(x^\alpha) $$

Now we calculate the Jumarie type fractional order derivative of $\exp(x) = E_1(x)$ like we did for $E_\alpha(x^\alpha)$ by using the formula ${}_0^J D_x^\alpha \left[ x^\upsilon \right] = \frac{\Gamma(1+\upsilon)}{\Gamma(1+\upsilon-\alpha)} x^{\upsilon-\alpha}$ and ${}_0^J D_x^\alpha [1] = 0$.

$$
\begin{aligned}
{}_0^J D_x^\alpha [\exp(ax)] &= {}_0^J D_x^\alpha [E_1(ax)] = {}_0^J D_x^\alpha \left( 1 + \frac{ax}{\Gamma(2)} + \frac{a^2 x^2}{\Gamma(3)} + \frac{a^3 x^3}{\Gamma(4)} + \ldots \right) \\
&= 0 + \frac{ax^{1-\alpha}}{\Gamma(2-\alpha)} + \frac{a^2 x^{2-\alpha}}{\Gamma(3-\alpha)} + \frac{a^3 x^{3-\alpha}}{\Gamma(4-\alpha)} + \ldots \\
&= a E_{1,2-\alpha}(ax)
\end{aligned}
$$

On the other hand the Jumarie type fractional order derivative of $\cos(ax)$ is following, as we did for $\cos_\alpha(x^\alpha)$ by using the formula ${}_0^J D_x^\alpha \left[ x^\upsilon \right] = \frac{\Gamma(1+\upsilon)}{\Gamma(1+\upsilon-\alpha)} x^{\upsilon-\alpha}$ and ${}_0^J D_x^\alpha [1] = 0$.

$$
\begin{aligned}
{}_0^J D_x^\alpha [\cos(ax)] &= {}_0^J D_x^\alpha \left[ 1 - \frac{a^2 x^2}{\Gamma(3)} + \frac{a^4 x^4}{\Gamma(5)} + \frac{a^6 x^6}{\Gamma(7)} + \ldots \right] \\
&= 0 - \frac{a^2 x^{2-\alpha}}{\Gamma(3-\alpha)} + \frac{a^4 x^{4-\alpha}}{\Gamma(5-\alpha)} - \frac{a^6 x^{6-\alpha}}{\Gamma(7-\alpha)} + \ldots \\
&= -ax^{1-\alpha} \sin_{1,2-\alpha}(ax)
\end{aligned}
$$

We obtain

$$ {}_0^J D_x^\alpha [\cos(ax)] = -ax^{1-\alpha} \sin_{1,2-\alpha}(ax) $$

Similarly the Jumarie type fractional order derivative of $\sin(x)$ is

$$_0^JD_x^\alpha[\sin(ax)] = ax^{1-\alpha}\cos_{1,2-\alpha}(ax)$$

## 2.1 Some definitions of roughness index

### a) Lipschitz Holder exponent(LHE)

A function is said to have LHE [1] as $\alpha$ if the following condition is satisfied

$$|f(x)-f(y)| \sim |x-y|^\alpha \qquad 0 < |x-y| < \varepsilon$$

Where $\varepsilon$ is a small positive number. The property LHE is a local property. The global LHE in interval $[a,b]$ is denoted by $\lambda$ and is defined by

$$\lambda = \inf_{x \in [a,b]} \alpha$$

unless $f(x)$ is a constant function, $\lambda \leq 1$. The Lipschitz Holder exponent is sometimes named as Holder exponent. For the continuous function $f: \mathbb{R} \to \mathbb{R}$, $f(x)$ satisfies the Lipschitz condition on its domain of definition $|f(x)-f(y)| < C|x-y|$ when $0 < |x-y| < \varepsilon$ where $\varepsilon$ is small positive number, and $C > 0$ is real constant. This function $f(x)$ has Holder exponent as $1$.

Consider the function

$f: \mathbb{R} \to \mathbb{R}$ such that $f(x) = \sin(x)$ then $|f(x)-f(y)| = |\sin(x)-\sin(y)| < C|x-y|$ when $0 < |x-y| < \varepsilon$ is a function with Holder exponent $1$. In a way it states that the continuous function in consideration is one-whole differentiable and the value of differentiation is bounded, that is $\frac{|f(x)-f(y)|}{|x-y|} < C$ for $0 < |x-y| < \varepsilon$.

### b) Holder Continuity

A continuous function $f(x)$ which is non-differentiable in classical sense is said to holder continuous with exponent $\alpha$ if

$$|f(x)-f(y)| < C|x-y|^\alpha \qquad 0 < |x-y| < \varepsilon$$

where $C > 0$ is a real constant and $\varepsilon > 0$.

### c) Fractional dimension

Fractional dimension ($d$) or box dimension [1] of a function or graph is local property, denotes the degree of roughness of a function or graph. Let the graph of a function is $f(x)$ for $x \in [a,b]$ can be covered by $N$-squares of size $r$ then with $\lim(r \to 0)$ the fractional dimension of the graph is defined as

$$d = \lim_{r \to 0} \frac{\log(N)}{\log(1/r)}$$

Again if $H$ be the Hurst exponent then the relation between the above Holder exponents are $\alpha = \lambda = H \quad d = 2 - H = 2 - \alpha$ [1], [9]. The Holder and Hurst exponents are equivalent for uni-fractal graphs that has a constant fractional dimension in defined interval [1], [9].

### d) Weierstrass function

In 1872 K. Weierstrass [23-25] proposed his famous example of an always continuous but nowhere differentiable function $W(x)$ on the real line $\mathbb{R}$ with two parameters $b \geq a > 1$ in the following form

$$W(x) = \sum_{k=1}^{\infty} a^{-k} \sin(b^k x) \qquad x \in \mathbb{R}$$

Where $b$ is odd-integer number. He proved that this function is continuous for all real number and is non-differentiable for all real values of $x$ if $ab > 1 + \frac{3\pi}{2}$. Considering $b$ a constant say $b = \lambda$, and $s = 2 - \frac{\log a}{\log b}$ another presentation of the Weierstrass function [13] can be obtained which is

$$W(x) = \sum_{k=1}^{\infty} \lambda^{(s-2)k} \sin(\lambda^k x) \qquad \lambda > 1 \qquad 1 < s < 2 \qquad (4)$$

In [13] Falconer established the fractional dimension of Weierstrass function defined in (4) is $s$ and the corresponding Holder exponent is $2 - s$.

### 3.0 The fractional Weierstrass Function

The original Weierstrass Function (4) is defined in the following form

$$W(x) = \sum_{k=1}^{\infty} \lambda^{(s-2)k} \sin(\lambda^k x) \qquad \lambda > 1 \qquad 1 < s < 2$$

We define the fractional Weierstrass Function in terms of Jumarie [2008] fractional sine function, that is $\sin_\alpha(x^\alpha)$ in the following form for $x \geq 0$

$$W_\alpha(x^\alpha) = \sum_{k=1}^{\infty} \lambda^{(s-2)k} \sin_\alpha(\lambda^{\alpha k} x^\alpha) \qquad \lambda > 1 \qquad 1 < s < 2$$

Where, $0 < \alpha < 1$, and for $\alpha = 1$ it becomes the original Weierstrass Function, and a condition that $W_\alpha(x^\alpha) = 0$ for $x < 0$.

We only are stating some lemmas which will be used to characterize the fractional Weierstrass function and its fractional derivative.

**Lemma 1:**

Let $f$ be function continuous on interval $[0,1]$ and $0 \leq s \leq 1$ [12-14],

Suppose

(1) $\qquad |f(x)-f(y)| \leq C|x-y|^s \qquad 0 < x \qquad y < 1$

then the dimension [12-14] of the graph $f$ is $d \leq 2-s$.

(2) Suppose $\delta_0 > 0$. For every $x \in [0,1]$, and $0 < \delta < \delta_0$ there exists $y \in [0,1]$ such that $|x-y| < \delta$ and $|f(x)-f(y)| \geq C\delta^s$ then the dimension [12-14] of the graph $f$ is $d \geq 2-s$.

**Theorem 1:** The Holder exponent of fractional Weierstrass function $W_\alpha(x^\alpha)$ with $0 < \alpha < 1$ is $2-s$ and consequently the Hausdorff dimension or fractional dimension is $s$ over any finite interval that is $[0,1]$.

**Proof:** We calculate $W_\alpha[(x+h)^\alpha] - W_\alpha[x^\alpha]$ in following steps where we have used our derived expression $\sin_\alpha\left(a(x+y)^\alpha\right) = \left[\sin_\alpha(ax^\alpha)\cos_\alpha(ay^\alpha) + \cos_\alpha(ax^\alpha)\sin_\alpha(ay^\alpha)\right]$

$$W_\alpha[(x+h)^\alpha] - W_\alpha[x^\alpha] = \sum_{k=1}^{\infty} \lambda^{(s-2)k} \sin_\alpha(\lambda^{\alpha k}(x+h)^\alpha) - \sum_{k=1}^{\infty} \lambda^{(s-2)k} \sin_\alpha(\lambda^{\alpha k}(x)^\alpha)$$

$$= \sum_{k=1}^{\infty} \lambda^{(s-2)k} \left[\sin_\alpha(\lambda^{\alpha k} x^\alpha)\cos_\alpha(\lambda^{\alpha k} h^\alpha) + \cos_\alpha(\lambda^{\alpha k} x^\alpha)\sin_\alpha(\lambda^{\alpha k} h^\alpha)\right]$$

$$- \sum_{k=1}^{\infty} \lambda^{(s-2)k} \sin_\alpha(\lambda^{\alpha k}(x)^\alpha)$$

$$= \sum_{k=1}^{\infty} \lambda^{(s-2)k} \left[\sin_\alpha(\lambda^{\alpha k} x^\alpha)\left(\cos_\alpha(\lambda^{\alpha k} h^\alpha) - 1\right) + \cos_\alpha(\lambda^{\alpha k} x^\alpha)\sin_\alpha(\lambda^{\alpha k} h^\alpha)\right]$$

From the series expansion of $\sin_\alpha(\lambda^{\alpha k} x^\alpha)$ and $\cos_\alpha(\lambda^{\alpha k} x^\alpha)$ and also from the figure-1 and 2, it is clear that for small $x$, $\sin_\alpha(\lambda^{\alpha k} x^\alpha) \approx \lambda^{\alpha k} x^\alpha$ and $\cos_\alpha(\lambda^{\alpha k} x^\alpha) \approx 1$ also both $\left|\sin_\alpha(\lambda^{\alpha k} x^\alpha)\right|$ and $\left|\cos_\alpha(\lambda^{\alpha k} x^\alpha)\right|$ is less than or equal to 1. Therefore, with above observation that is for small $h$, $\sin_\alpha(\lambda^{\alpha k} h^\alpha) \approx \lambda^{\alpha k} h^\alpha$, $\cos_\alpha(\lambda^{\alpha k} h^\alpha) - 1 \approx 0$ and for large $h$, $\cos_\alpha(\lambda^{\alpha k} h^\alpha) \approx 0$ we write the following

$$|W_\alpha[(x+h)^\alpha]-W_\alpha[x^\alpha]| \le \sum_{k=1}^{\infty} \lambda^{(s-2)k}\Big[|\sin_\alpha(\lambda^{\alpha k}x^\alpha)||\cos_\alpha(\lambda^{\alpha k}h^\alpha)-1|+|\cos_\alpha(\lambda^{\alpha k}x^\alpha)||\sin_\alpha(\lambda^{\alpha k}h^\alpha)|\Big]$$

$$\le \sum_{k=1}^{\infty} \lambda^{(s-2)k}\Big[\min(\lambda^{\alpha k}h^\alpha,1)\Big]$$

Choose $0 < h < 1$ then one can find positive integer $m$ such that $\lambda^{-(m+1)} \le h \le \lambda^{-m}$ then divide the summation that is $\sum_{k=1}^{\infty}\lambda^{(s-2)k}\Big[\min(\lambda^{\alpha k}h^\alpha,1)\Big]$ into two parts. First part for $k=1$ to $m$ then $\sin_\alpha(\lambda^{\alpha k}x^\alpha) \approx \lambda^{\alpha k}x^\alpha$ and for other values of $k$ maximum value of the expression in third bracket is equal to 1. We use the geometric series formulas $\sum_{k=1}^{m}a^k = a\left(\frac{a^m-1}{a-1}\right)$ and $\sum_{k=1}^{\infty}a^k = \frac{a}{1-a}$, for $\sum_{k=m+1}^{\infty}x^k = \frac{x^{m+1}}{1-x}$ in the following derivation.

$$|W_\alpha[(x+h)^\alpha]-W_\alpha[x^\alpha]| \le \sum_{k=1}^{m}\lambda^{(s-2)k}\left(\lambda^{\alpha k}h^\alpha\right)+1\sum_{k=m+1}^{\infty}\lambda^{(s-2)k}$$

$$= h^\alpha \sum_{k=1}^{m}\lambda^{(s-2+\alpha)k} + \sum_{k=m+1}^{\infty}\lambda^{(s-2)k}$$

$$= h^\alpha\left(\lambda^{(s-2+\alpha)}\frac{\lambda^{(s-2+\alpha)m}-1}{\lambda^{(s-2+\alpha)}-1}\right)+\left(\frac{\lambda^{(s-2)(m+1)}}{1-\lambda^{(s-2)}}\right)$$

$$\le h^\alpha \frac{\lambda^{(s-2+\alpha)(m+1)}}{\lambda^{(s-2+\alpha)}-1}+\frac{\lambda^{(s-2)(m+1)}}{1-\lambda^{(s-2)}}$$

With $\lambda^{-(m+1)} \le h \le \lambda^{-m}$, that is $\lambda^{(m+1)} \ge h^{-1} \ge \lambda^m$ we get the following

$$|W_\alpha[(x+h)^\alpha]-W_\alpha[x^\alpha]| \le h^\alpha \frac{h^{-(s-2+\alpha)}}{\lambda^{(s-2+\alpha)}-1}+\frac{h^{-(s-2)}}{1-\lambda^{(s-2)}}$$

$$= \left(\frac{1}{\lambda^{(s-2+\alpha)}-1}+\frac{1}{1-\lambda^{(s-2)}}\right)h^{2-s} = C_1 h^{2-s}$$

Where the constant $C_1 = \frac{1}{\lambda^{(s-2+\alpha)}-1}+\frac{1}{1-\lambda^{(s-2)}}$. From definition of Holderian function and the above discussion it is clear that fractional Weierstrass function is also Holder continuous with Holder exponent $(2-s)$, a fractional number. This shows (by Lemma-1) that Hausdorff dimension of graph of fractional Weierstrass function is $[2-(2-s)]=s$. Thus the Hausdorff dimension of fractional Weierstrass function and original Weierstrass function is same, is independent of fractional exponent ($\alpha$) as defined in (4).

## 4.0 The Jumarie fractional derivative of fractional Weierstrass Function

Many authors found the fractional derivative of the continuous but nowhere differentiable function that is Weierstrass Function [10-17] using different type definitions of fractional derivatives. Here Jumarie type fractional order derivative of $W_\alpha(x^\alpha)$ is of order $\alpha$

$$\,_0^J D_x^\alpha \left[ W_\alpha(x^\alpha) \right] = \sum_{k=1}^{\infty} \lambda^{(s-2)k} \left( \,_0^J D_x^\alpha \left[ \sin_\alpha(\lambda^{\alpha k} x^\alpha) \right] \right)$$

$$= \sum_{k=1}^{\infty} \lambda^{(s-2)k} \lambda^{\alpha k} \cos_\alpha(\lambda^{\alpha k} x^\alpha)$$

We used in above derivation the identity $\,_0^J D_x^\alpha \left[ \sin_\alpha(ax^\alpha) \right] = a \cos_\alpha(ax^\alpha)$. Therefore from above derivation we obtain the following

$$\,_0^J D_x^\alpha \left[ W_\alpha(x^\alpha) \right] = \sum_{k=1}^{\infty} \lambda^{(s-2+\alpha)k} \cos_\alpha(\lambda^{\alpha k} x^\alpha) \qquad (5)$$

Since if $0 < \alpha < 1$ then $\cos_\alpha(\lambda^{\alpha k} x^\alpha)$ is a bounded function and therefore $\,_0^J D_x^\alpha \left[ W_\alpha(x^\alpha) \right]$ will be bounded function if $\sum_{k=1}^{\infty} \lambda^{(s-2+\alpha)k}$ is convergent. Since $\sum_{k=1}^{\infty} \lambda^{(s-2+\alpha)k}$ is a geometric series will be convergent if $s - 2 + \alpha < 0$ implying $\alpha < 2 - s$. Hence the fractional derivative of order $\alpha$ with $0 < \alpha < 1$ of the Weierstrass Function will exists when $\alpha < 2 - s$.

Again if $\alpha > 1$ then $\cos_\alpha(\lambda^{\alpha k} x^\alpha)$ and $\sin_\alpha(\lambda^{\alpha k} x^\alpha)$ for $k = 1, 2, 3...$; are unbounded functions (figure-1 and 2) and will grow by oscillating without bound to $\pm\infty$ for $x \to \infty$. Since $1 < s < 2$ and $\alpha > 1$ implying $s + \alpha - 2 > 0$ therefore $\sum_{k=1}^{\infty} \lambda^{(s-2+\alpha)k}$ is a divergent series. Therefore

$$\,_0^J D_x^\alpha \left[ W_\alpha(x^\alpha) \right] = \sum_{k=1}^{\infty} \lambda^{(s-2+\alpha)k} \cos_\alpha(\lambda^{\alpha k} x^\alpha)$$

is a divergent series for $\alpha > 1$. We write following observation

$$\,_0^J D_x^\alpha \left[ W_\alpha(x^\alpha) \right] = \begin{cases} \text{Bounded for} & \alpha < 2 - s \\ \text{Unbounded for} & \alpha \geq 2 - s \end{cases}$$

This shows that $\alpha$-order ($0 < \alpha < 1$) fractional derivative of the fractional Weierstrass function exists when $\alpha < 2 - s$ and for $\alpha \geq 2 - s$ it does not exist. Thus we can state a theorem in the following form

**Theorem 2:** $\alpha$-order ($0 < \alpha < 1$) fractional derivative of the fractional Weierstrass function

$$W_\alpha(x^\alpha) = \sum_{k=1}^{\infty} \lambda^{(s-2)k} \sin_\alpha(\lambda^{\alpha k} x^\alpha) \qquad \lambda > 1 \qquad 1 < s < 2$$

exists when $\alpha < 2 - s$ and for $\alpha \geq 2 - s$ it does not exist.

**Theorem 3:** The Holder exponent of $\alpha$-order fractional derivative of fractional Weierstrass function $W_\alpha(x^\alpha)$, $0 < \alpha < 1$ is $2 - s - \alpha$ and consequently the Hausdorff dimension or fractional dimension is $s + \alpha$ over any finite interval $[0,1]$.

**Proof:** Let

$$_0^J D_x^\alpha \left[ W_\alpha(x^\alpha) \right] = W_\alpha^{(\alpha)}(x^\alpha) = \sum_{k=1}^{\infty} \lambda^{(s-2+\alpha)k} \cos_\alpha(\lambda^{\alpha k} x^\alpha)$$

denotes $\alpha$-order fractional Jumarie derivative of fractional Weierstrass function, Then using the identity $\cos_\alpha\left(a(x+y)^\alpha\right) = \left[\cos_\alpha(ax^\alpha)\cos_\alpha(ay^\alpha) - \sin_\alpha(ax^\alpha)\sin_\alpha(ay^\alpha)\right]$ we get the following

$$W_\alpha^{(\alpha)}[(x+h)^\alpha] - W_\alpha^{(\alpha)}[x^\alpha] = \sum_{k=1}^{\infty} \lambda^{(s-2+\alpha)k} \cos_\alpha(\lambda^{\alpha k}(x+h)^\alpha) - \sum_{k=1}^{\infty} \lambda^{(s-2+\alpha)k} \cos_\alpha(\lambda^{\alpha k}(x)^\alpha)$$

$$= \sum_{k=1}^{\infty} \lambda^{(s-2+\alpha)k} \left[ \cos_\alpha(\lambda^{\alpha k} x^\alpha)\cos_\alpha(\lambda^{\alpha k} h^\alpha) - \sin_\alpha(\lambda^{\alpha k} x^\alpha)\sin_\alpha(\lambda^{\alpha k} h^\alpha) \right]$$

$$- \sum_{k=1}^{\infty} \lambda^{(s-2+\alpha)k} \cos_\alpha(\lambda^{\alpha k}(x)^\alpha)$$

$$= \sum_{k=1}^{\infty} \lambda^{(s-2+\alpha)k} \left[ \cos_\alpha(\lambda^{\alpha k} xx^\alpha)\left(\cos_\alpha(\lambda^{\alpha k} h^\alpha) - 1\right) - \sin_\alpha(\lambda^{\alpha k} x^\alpha)\sin_\alpha(\lambda^{\alpha k} h^\alpha) \right]$$

From the series expansion of $\sin_\alpha(\lambda^{\alpha k} x^\alpha)$ and $\cos_\alpha(\lambda^{\alpha k} x^\alpha)$ and also from the figure-1 and 2 it is clear that for small $x$, $\sin_\alpha(\lambda^{\alpha k} x^\alpha) \approx \lambda^{\alpha k} x^\alpha$ and $\cos_\alpha(\lambda^{\alpha k} x^\alpha) \approx 1$ also both $\left|\sin_\alpha(\lambda^{\alpha k} x^\alpha)\right|$ and $\left|\cos_\alpha(\lambda^{\alpha k} x^\alpha)\right|$ is less than or equal to 1. Therefore, with above observation that is for small $h$, $\sin_\alpha(\lambda^{\alpha k} h^\alpha) \approx \lambda^{\alpha k} h^\alpha$, $\cos_\alpha(\lambda^{\alpha k} h^\alpha) - 1 \approx 0$ and for large $h$, $\cos_\alpha(\lambda^{\alpha k} h^\alpha) \approx 0$ we write the following

$$|W_\alpha^{(\alpha)}[(x+h)^\alpha] - W_\alpha^{(\alpha)}[x^\alpha]| \leq \sum_{k=1}^{\infty} \lambda^{(s-2+\alpha)k} [|\cos_\alpha(\lambda^{\alpha k} x^\alpha)||\cos_\alpha(\lambda^{\alpha k} h^\alpha) - 1| + |\sin_\alpha(\lambda^{\alpha k} x^\alpha)||\sin_\alpha(\lambda^{\alpha k} h^\alpha)|]$$

$$\leq \sum_{k=1}^{\infty} \lambda^{(s-2+\alpha)k} \left[ \min\left(\lambda^{\alpha k} h^\alpha, 1\right) \right]$$

Choose $0 < h < 1$ then one can find positive integer $m$ such that $\lambda^{-(m+1)} \leq h \leq \lambda^{-m}$ then as per our earlier derivation for $W_\alpha(x^\alpha)$ we do the following steps

$$|W_\alpha^{(\alpha)}[(x+h)^\alpha] - W_\alpha^{(\alpha)}[x^\alpha]| \leq \sum_{k=1}^{m} \lambda^{(s-2+\alpha)k} \lambda^{\alpha k} h^\alpha + \sum_{k=m+1}^{\infty} (1)\lambda^{(s-2+\alpha)k}$$

$$= h^\alpha \sum_{k=1}^{m} \lambda^{(s-2+2\alpha)} + \sum_{k=m+1}^{\infty} (1)\lambda^{(s-2+\alpha)k}$$

$$= h^\alpha \left( \lambda^{(s-2+2\alpha)} \frac{\lambda^{(s-2+2\alpha)m} - 1}{\lambda^{(s-2+2\alpha)} - 1} \right) + 1 \left( \frac{\lambda^{(s-2+\alpha)(m+1)}}{1 - \lambda^{(s-2+\alpha)}} \right)$$

$$\leq h^\alpha \frac{\lambda^{(s-2+2\alpha)(m+1)}}{\lambda^{(s-2+2\alpha)} - 1} + \frac{\lambda^{(s-2+\alpha)(m+1)}}{1 - \lambda^{(s-2+\alpha)}}$$

With $\lambda^{-(m+1)} \leq h \leq \lambda^{-m}$, that is $\lambda^{(m+1)} \geq h^{-1} \geq \lambda^m$ we get the following

$$\left| W_\alpha^{(\alpha)}[(x+h)^\alpha] - W_\alpha^{(\alpha)}[x^\alpha] \right| \leq h^\alpha \frac{\lambda^{(s-2+2\alpha)(m+1)}}{\lambda^{(s-2+2\alpha)} - 1} + \frac{\lambda^{(s-2+\alpha)(m+1)}}{1 - \lambda^{(s-2+\alpha)}}$$

$$\leq h^\alpha \frac{h^{-(s-2+2\alpha)}}{\lambda^{(s-2+2\alpha)} - 1} + \frac{h^{-(s-2+\alpha)}}{1 - \lambda^{(s-2+\alpha)}}$$

$$\leq \left( \frac{1}{\lambda^{(s-2+2\alpha)} - 1} + \frac{1}{1 - \lambda^{(s-2+\alpha)}} \right) h^{2-s-\alpha}$$

$$\leq C_2 h^{2-s-\alpha}$$

Where $C_2 = \frac{1}{\lambda^{(s-2+2\alpha)} - 1} + \frac{1}{1 - \lambda^{(s-2+\alpha)}}$. From definition of Holderian function and above discussion it is clear that $\alpha$-order ($0 < \alpha < 1$) fractional derivative of fractional Weierstrass function is also Holder continuous with Holder exponent $2 - s - \alpha$. This shows that Hausdorff dimension of graph of fractional Weierstrass function is $[2 - (2 - s - \alpha)] = s + \alpha$ (by lemma-1). The graph dimension increased by fractional order for fractional derivative of Weierstrass function by amount of fractional derivative-the graph becomes rougher.

## Conclusion

The fractional Weierstrass function is a continuous function for all real values of the arguments and its box dimension and Holder exponent is independent of fractional order that incorporated to the fractional Weierstrass functions. Again the Box dimension of fractional derivative of the fractional Weierstrass increases with increase of order of fractional derivative. This invariant nature of the roughness index of fractional Weierstrass function when generalized with fractional trigonometric function is remarkable. The other embodiment in similar lines as in this paper to get different fractional Weierstrass function is under development.

## Acknowledgement

Acknowledgments are to **Board of Research in Nuclear Science** (BRNS), Department of Atomic Energy Government of India for financial assistance received through BRNS research project no. 37(3)/14/46/2014-BRNS with BSC BRNS, title "Characterization of unreachable (Holderian) functions via Local Fractional Derivative and Deviation Function".